\documentclass{gtart}

\def\ifplaintex{\expandafter\ifx\csname documentclass\endcsname\relax}

\def\gtp{{\mathsurround=0pt\it $\cal G\mskip-2mu$eometry \&\ 
$\cal T\!\!$opology $\cal P\!$ublications}}  

\def\recd{{\small Received:\qua\receiveddate\ifx\reviseddate\relax
\else\qquad Revised:\qua\reviseddate\fi\par}} 

\def\lognumber#1{\def\thelognumber{#1}}
\def\volumenumber#1{\def\thevolumenumber{#1}}
\def\volumeyear#1{\def\thevolumeyear{#1}}
\def\papernumber#1{\def\thepapernumber{#1}}
\def\pagenumbers#1#2{\def\startpage{#1}\def\finishpage{#2}}
\def\published#1{\def\publishdate{#1}}

\def\received#1{\def\receiveddate{#1}}
\def\revised#1{\def\reviseddate{#1}}
\def\accepted#1{\def\accepteddate{#1}}

\long\def\asciiabstract#1{\long\def\theasciiabstract{#1}}

\let\\\par\let\thelognumber\relax\let\thevolumenumber\relax
\let\thepapernumber\relax\let\thevolumeyear\relax\let\startpage\relax
\let\finishpage\relax\let\publishdate\relax\let\receiveddate\relax
\let\reviseddate\relax\let\accepteddate\relax\let\theasciititle\relax
\let\theasciiauthors\relax
\let\theasciiabstract\relax

\let\theasciiemail\relax

\ifplaintex
\font\logobig=cmssbx10 scaled 3836
\font\logomed=cmssbx10 scaled 2557
\else
\font\logobig=cmssbx10 scaled 4200
\font\logomed=cmssbx10 scaled 2800
\fi

\long\def\makeagttitle{   
\count0=\startpage
\agt\hfill      
\hbox to 45truept{\vbox to 0pt{\vglue -13truept{\logomed A\kern -.37em{\logobig 
T}\kern -.38em G}\vss}\hss}
\break
{\small Volume \thevolumenumber\ (\thevolumeyear)
\startpage--\finishpage\nl
Published: \publishdate}

\vglue .25truein

{\parskip=0pt\leftskip 0pt plus
1fil\def\\{\par\smallskip}{\Large\bf\thetitle}\par\medskip} \vglue
0.05truein

{\parskip=0pt\leftskip 0pt plus 1fil\def\\{\par}{\sc\theauthors}
\par\medskip}%
 
\vglue 0.03truein 

{\small\leftskip 25truept\rightskip 25truept{\bf Abstract}\stdspace\theabstract

{\bf AMS Classification}\stdspace\theprimaryclass
\ifx\thesecondaryclass\relax\else; \thesecondaryclass\fi\par
{\bf Keywords}\stdspace \thekeywords\par}\vglue 7truept

}   

\ifplaintex
\font\phead=cmsl9 scaled 950
\font\pnum=cmbx10 scaled 913
\font\pfoot=cmsl9 scaled 950
\headline{\vbox to 0pt{\vskip -4.5mm\line{\small\phead\ifnum
\count0=\startpage ISSN 1472-2739 (on-line) 1472-2747 (printed)
\hfill {\pnum\folio}\else\ifodd\count0\def\\{ }%
\ifx\theshorttitle\relax\thetitle\else\theshorttitle\fi\hfill{\pnum\folio}
\else\def\\{ and }{\pnum\folio}\hfill\ifx\theshortauthors\relax\theauthors
\else\theshortauthors\fi\fi\fi}\vss}}
\footline{\vbox to 0pt{\vglue 0mm\line{\small\pfoot\ifnum\count0=\startpage
\copyright\ \gtp\hfill\else
\agt, Volume \thevolumenumber\ (\thevolumeyear)\hfill\fi}\vss}}
\else
\font=cmsl9 scaled 1050
\font\lnum=cmbx10 
\font=cmsl9 scaled 1050
\makeatletter
\def\@oddhead{{\small\lhead\ifnum\count0=\startpage ISSN 1472-2739 
(on-line) 1472-2747 (printed)\hfill {\lnum\number\count0}\else\ifodd\count0
\def\\{ }\ifx\theshorttitle\relax \thetitle \else\theshorttitle\fi\hfill
{\lnum\number\count0}\else\def\\{ and }{\lnum\number\count0}
\hfill\ifx\theshortauthors\relax 
\theauthors\else\theshortauthors\fi\fi\fi}}\def\@evenhead{\@oddhead}
\def\@oddfoot{\small\lfoot\ifnum\count0=\startpage\copyright\ \gtp\hfill\else
\agt, Volume \thevolumenumber\ (\thevolumeyear)\hfill\fi}
\def\@evenfoot{\@oddfoot}
\makeatother
\fi
\let\maketitlepage\makeagttitle
\let\makeshorttitle\maketitlepage
\let\maketitle\maketitlepage

\newwrite\gtoutfile
\long\gdef\makeheadfile{  
{\def\\{, }\def\s{ }
\immediate\openout\gtoutfile head.xxx
\immediate\write\gtoutfile{To: math@arxiv.org}
\immediate\write\gtoutfile{Subject: put OR rep NNNNN:ppppp}
\immediate\write\gtoutfile{--text follows this line--}
\immediate\write\gtoutfile{Proxy-for: \ifx\theasciiauthors\relax
\theauthors\else\theasciiauthors\fi\s<\ifx\theasciiemail\relax\theemail\else\theasciiemail\fi>}
\immediate\write\gtoutfile{\noexpand\\}
\immediate\write\gtoutfile{Authors: \ifx\theasciiauthors\relax
\theauthors\else\theasciiauthors\fi}
{\def\\{ }\immediate\write\gtoutfile{Title: \ifx\theasciititle\relax
\thetitle\else\theasciititle\fi}}
\immediate\write\gtoutfile{Subj-class: GT or SG, GR etc}
\immediate\write\gtoutfile{MSC-class: \theprimaryclass\ifx\thesecondaryclass\relax\else, \thesecondaryclass\fi}
\immediate\write\gtoutfile{Journal-ref: Algebr. Geom. Topol. \thevolumenumber\s
(\thevolumeyear) \startpage-\finishpage}
\immediate\write\gtoutfile{Comments: Published by Algebraic and
Geometric Topology at}
\immediate\write\gtoutfile{\s\s\s  http://www.maths.warwick.ac.uk/agt/AGTVol\thevolumenumber/agt-\thevolumenumber-\thepapernumber.abs.html}
\immediate\write\gtoutfile{\noexpand\\}
\immediate\write\gtoutfile{}
\ifx\theasciiabstract\relax
\immediate\write\gtoutfile{\theabstract}\else
\immediate\write\gtoutfile{\theasciiabstract}\fi
\immediate\write\gtoutfile{}
\immediate\write\gtoutfile{\noexpand\\}
\immediate\write\gtoutfile{}
\immediate\closeout\gtoutfile}}  

\def\maketitlepage{\makeagttitle\makeheadfile}
\let\makeshorttitle\maketitlepage
\let\maketitle\maketitlepage

\def\ifplaintex{\expandafter\ifx\csname documentclass\endcsname\relax}

\def\gtp{{\mathsurround=0pt\it $\cal G\mskip-2mu$eometry \&\ 
$\cal T\!\!$opology $\cal P\!$ublications}}  

\def\recd{{\small Received:\qua\receiveddate\ifx\reviseddate\relax
\else\qquad Revised:\qua\reviseddate\fi\par}} 

\def\lognumber#1{\def\thelognumber{#1}}
\def\volumenumber#1{\def\thevolumenumber{#1}}
\def\volumeyear#1{\def\thevolumeyear{#1}}
\def\papernumber#1{\def\thepapernumber{#1}}
\def\pagenumbers#1#2{\def\startpage{#1}\def\finishpage{#2}}
\def\published#1{\def\publishdate{#1}}

\def\received#1{\def\receiveddate{#1}}
\def\revised#1{\def\reviseddate{#1}}
\def\accepted#1{\def\accepteddate{#1}}

\long\def\asciiabstract#1{\long\def\theasciiabstract{#1}}

\let\\\par\let\thelognumber\relax\let\thevolumenumber\relax
\let\thepapernumber\relax\let\thevolumeyear\relax\let\startpage\relax
\let\finishpage\relax\let\publishdate\relax\let\receiveddate\relax
\let\reviseddate\relax\let\accepteddate\relax\let\theasciititle\relax
\let\theasciiauthors\relax
\let\theasciiabstract\relax

\let\theasciiemail\relax

\ifplaintex
\font\logobig=cmssbx10 scaled 3836
\font\logomed=cmssbx10 scaled 2557
\else
\font\logobig=cmssbx10 scaled 4200
\font\logomed=cmssbx10 scaled 2800
\fi

\long\def\makeagttitle{   
\count0=\startpage
\agt\hfill      
\hbox to 45truept{\vbox to 0pt{\vglue -13truept{\logomed A\kern -.37em{\logobig 
T}\kern -.38em G}\vss}\hss}
\break
{\small Volume \thevolumenumber\ (\thevolumeyear)
\startpage--\finishpage\nl
Published: \publishdate}

\vglue .25truein

{\parskip=0pt\leftskip 0pt plus
1fil\def\\{\par\smallskip}{\Large\bf\thetitle}\par\medskip} \vglue
0.05truein

{\parskip=0pt\leftskip 0pt plus 1fil\def\\{\par}{\sc\theauthors}
\par\medskip}%
 
\vglue 0.03truein 

{\small\leftskip 25truept\rightskip 25truept{\bf Abstract}\stdspace\theabstract

{\bf AMS Classification}\stdspace\theprimaryclass
\ifx\thesecondaryclass\relax\else; \thesecondaryclass\fi\par
{\bf Keywords}\stdspace \thekeywords\par}\vglue 7truept

}   

\ifplaintex
\font\phead=cmsl9 scaled 950
\font\pnum=cmbx10 scaled 913
\font\pfoot=cmsl9 scaled 950
\headline{\vbox to 0pt{\vskip -4.5mm\line{\small\phead\ifnum
\count0=\startpage ISSN 1472-2739 (on-line) 1472-2747 (printed)
\hfill {\pnum\folio}\else\ifodd\count0\def\\{ }%
\ifx\theshorttitle\relax\thetitle\else\theshorttitle\fi\hfill{\pnum\folio}
\else\def\\{ and }{\pnum\folio}\hfill\ifx\theshortauthors\relax\theauthors
\else\theshortauthors\fi\fi\fi}\vss}}
\footline{\vbox to 0pt{\vglue 0mm\line{\small\pfoot\ifnum\count0=\startpage
\copyright\ \gtp\hfill\else
\agt, Volume \thevolumenumber\ (\thevolumeyear)\hfill\fi}\vss}}
\else
\font=cmsl9 scaled 1050
\font\lnum=cmbx10 
\font=cmsl9 scaled 1050
\makeatletter
\def\@oddhead{{\small\lhead\ifnum\count0=\startpage ISSN 1472-2739 
(on-line) 1472-2747 (printed)\hfill {\lnum\number\count0}\else\ifodd\count0
\def\\{ }\ifx\theshorttitle\relax \thetitle \else\theshorttitle\fi\hfill
{\lnum\number\count0}\else\def\\{ and }{\lnum\number\count0}
\hfill\ifx\theshortauthors\relax 
\theauthors\else\theshortauthors\fi\fi\fi}}\def\@evenhead{\@oddhead}
\def\@oddfoot{\small\lfoot\ifnum\count0=\startpage\copyright\ \gtp\hfill\else
\agt, Volume \thevolumenumber\ (\thevolumeyear)\hfill\fi}
\def\@evenfoot{\@oddfoot}
\makeatother
\fi
\let\maketitlepage\makeagttitle
\let\makeshorttitle\maketitlepage
\let\maketitle\maketitlepage

\newwrite\gtoutfile
\long\gdef\makeheadfile{  
{\def\\{, }\def\s{ }
\immediate\openout\gtoutfile head.xxx
\immediate\write\gtoutfile{To: math@arxiv.org}
\immediate\write\gtoutfile{Subject: put OR rep NNNNN:ppppp}
\immediate\write\gtoutfile{--text follows this line--}
\immediate\write\gtoutfile{Proxy-for: \ifx\theasciiauthors\relax
\theauthors\else\theasciiauthors\fi\s<\ifx\theasciiemail\relax\theemail\else\theasciiemail\fi>}
\immediate\write\gtoutfile{\noexpand\\}
\immediate\write\gtoutfile{Authors: \ifx\theasciiauthors\relax
\theauthors\else\theasciiauthors\fi}
{\def\\{ }\immediate\write\gtoutfile{Title: \ifx\theasciititle\relax
\thetitle\else\theasciititle\fi}}
\immediate\write\gtoutfile{Subj-class: GT or SG, GR etc}
\immediate\write\gtoutfile{MSC-class: \theprimaryclass\ifx\thesecondaryclass\relax\else, \thesecondaryclass\fi}
\immediate\write\gtoutfile{Journal-ref: Algebr. Geom. Topol. \thevolumenumber\s
(\thevolumeyear) \startpage-\finishpage}
\immediate\write\gtoutfile{Comments: Published by Algebraic and
Geometric Topology at}
\immediate\write\gtoutfile{\s\s\s  http://www.maths.warwick.ac.uk/agt/AGTVol\thevolumenumber/agt-\thevolumenumber-\thepapernumber.abs.html}
\immediate\write\gtoutfile{\noexpand\\}
\immediate\write\gtoutfile{}
\ifx\theasciiabstract\relax
\immediate\write\gtoutfile{\theabstract}\else
\immediate\write\gtoutfile{\theasciiabstract}\fi
\immediate\write\gtoutfile{}
\immediate\write\gtoutfile{\noexpand\\}
\immediate\write\gtoutfile{}
\immediate\closeout\gtoutfile}}  

\def\maketitlepage{\makeagttitle\makeheadfile}
\let\makeshorttitle\maketitlepage
\let\maketitle\maketitlepage

\lognumber{13}
\volumenumber{1}
\volumeyear{2001}
\papernumber{13}
\pagenumbers{271}{297}
\received{12 February 2001}
\revised{10 May 2001}
\accepted{16 May 2001}
\published{22 May 2001}

\usepackage{amsmath}
\usepackage{amssymb}
\usepackage{latexsym}
\usepackage[all]{xy}

\newcommand{\noproof}{\qed}

\theoremstyle{plain}
\newtheorem{theorem}{Theorem}[section]
\newtheorem{proposition}[theorem]{Proposition}
\newtheorem{lemma}[theorem]{Lemma}
\newtheorem{corollary}[theorem]{Corollary}
\theoremstyle{definition}
\newtheorem{definition}[theorem]{Definition}

\newcommand{\supp}{\mathop{supp}}

\begin{document}

\title{Coarse homology theories}

\author{Paul D. Mitchener}

\address{Department of Mathematics and Computer Science, Odense University,
Denmark}

\email{mitchene@imada.sdu.dk}

\url{http://www.imada.sdu.dk/\char'176mitchene/index2.html}

\begin{abstract}
In this paper we develop an axiomatic approach to coarse homology theories.  We prove a uniqueness result concerning coarse homology theories on the category of ``coarse $CW$-complexes''.  This uniqueness result is used to prove a version of the coarse Baum-Connes conjecture for such spaces.
\end{abstract}

\asciiabstract{
In this paper we develop an axiomatic approach to coarse homology theories.  We prove a uniqueness result concerning coarse homology theories on the category of `coarse CW-complexes'.  This uniqueness result is used to prove a version of the coarse Baum-Connes conjecture for such spaces.}

\primaryclass{55N35, 55N40}
\secondaryclass{19K56, 46L85}
\keywords{Coarse geometry, exotic homology, coarse Baum-Connes conjecture, Novikov conjecture}

\maketitle

\section{Introduction}

Coarse geometry is the study of large scale structures on geometric objects.  A large-scale analogue of cohomology, called coarse cohomology, was introduced by J.~Roe in \cite{Roe3} in order to perform index theory on non-compact manifolds.  The dual theory, coarse homology, has proved useful in formulating a coarse version of the Baum-Connes conjecture; see \cite{HR2} and \cite{Yu}.\footnote{Strictly speaking, coarse $K$-homology rather than the coarse version of ordinary homology is used here.}  This {\em coarse Baum-Connes conjecture} has similar consequences to the original Baum-Connes conjecture formulated in \cite{BCH}.  In particular, the coarse Baum-Connes conjecture implies the Novikov conjecture concerning the oriented homotopy-invariance of higher signatures on a manifold.  The book \cite{Roe1} contains an overview of coarse geometric techniques applied to such problems.

In this paper we present an axiomatic approach to coarse (generalised) homology theories; the axioms used are coarse analogues of the classical Eilenberg-Steenrod axioms (see \cite{ES}).

A {\em coarse $(n-1)$-sphere} is defined to be Euclidean space ${\mathbb R}^n$ equipped with a suitable coarse structure compatible with the topology.  This terminology comes from thinking of a coarse $(n-1)$-sphere as an ordinary $(n-1)$-sphere `at infinity'.  There is a similar notion of a {\em coarse $n$-cell}.  Using these basic building blocks we can formulate the notion of a {\em coarse $CW$-complex}.

The main result of this article is the fact that a coarse homology theory is uniquely determined on the category of finite coarse $CW$-complexes by its values on the one-point coarse space and on coarse $0$-cells.  The proof of this result is analogous to the corresponding uniqueness result for generalised homology theories in classical algebraic topology.

As an application, we use this result together with the main results from \cite{HPR} to study the coarse Baum-Connes conjecture.  In particular, using coarse structures arising from `continuous control at infinity' (see \cite{ACFP}) we obtain a result on the Novikov conjecture that appears to be new.

The remainder of the paper comprises: Section 2: The coarse category;
section 3: Coarse homology theories; section 4: Relative coarse
homology; section 5: Coarse $CW$-complexes; section 6: The coarse
assembly map.

\subsection*{Acknowledgements}

This article was written with support from TMR network ERB FMRX CT-97-0107 on algebraic $K$-theory, linear algebraic groups, and related structures.  Discussions with Hans J{\o}rgen Munkholm led to the concept of a coarse $CW$-complex and the generalised notion of coarsening introduced here.  Earlier discussions with John Roe indicated that the coarse Baum-Connes conjecture could be proved by showing that the functors $X\mapsto K_n C^\ast X$ form a coarse homology theory in some suitable sense.  More recently, he also alerted me to the article \cite{STY}.  Naturally, I am extremely grateful to both of these individuals.

\section{The coarse category}

The basic concepts of coarse geometry can be defined axiomatically; the following definition is a special case of the definition of a bornology on a topological space, $X$, given in \cite{HPR}.

\begin{definition}
A set $X$ is called a {\em coarse space} if there is a distinguished collection, $\cal E$, of subsets of the product $X\times X$ called {\em entourages} such that:

\begin{itemize}

\item Any finite union of entourages is contained in an entourage.

\item The union of all entourages is the entire space $X\times X$.

\item The {\em inverse} of an entourage $M$:
$$M^{-1} = \{ (y,x)\in X\times X \ |\ (x,y)\in M \}$$
is contained in an entourage.

\item The {\em composition} of entourages $M_1$ and $M_2$:
$$M_1 M_2 = \{ (x,z)\in X\times X \ |\ (x,y)\in M_1 \ (y,z)\in M_2 \textrm{ for some }y\in X \}$$
is contained in an entourage.

\end{itemize}

The coarse space $X$ is called {\em unital} if the diagonal, $D = \{ (x,x) \ |\ x\in X \}$, is contained in an entourage.

\end{definition}

Two coarse structures $\{ M_\alpha \subseteq X\times X \ |\ \alpha \in A \}$ and $\{ N_\beta \subseteq X\times X \ |\ \beta \in B \}$ on a set $X$ are said to be {\em equivalent} if every entourage $M_\alpha$ is contained in some entourage $N_\beta$, and every entourage $N_{\beta'}$ is contained in some entourage $M_{\alpha'}$.  We do not distinguish between equivalent coarse structures.

\begin{definition}
Let $X$ be a coarse space, and consider maps $f\co S\rightarrow X$ and $g\co S\rightarrow X$ for some set $S$.  Then the maps $f$ and $g$ are called {\em close} or {\em coarsely equivalent} if the set:
$$\{ (f(s) , g(s)) \ |\ s\in S \}$$
is contained in an entourage.
\end{definition}

There is a notion of a subset of a coarse space being of finite size.

\begin{definition} \label{bounded}
Let $X$ be a coarse space, and consider a subset $A\subseteq X$ and an entourage $M\subseteq X\times X$.  Then we define:
$$M[A] = \{ x\in X \ |\ (x,a)\in M \textrm{ for some point }a\in A \}$$

As a special case, for a single point $x\in X$ we write $M(x) = M[\{ x \}]$.  We call a subset $B\subseteq X$ {\em bounded} if it is contained in a set of the form $M(x)$.
\end{definition}

Observe that a subset of a bounded set is bounded, and a finite union of bounded sets is bounded.  If $B\subseteq X$ is a bounded set, and $M$ is an entourage, then the set $M[B]$ is bounded.

In this article we explore invariants that depend only on the coarse structure of a given space.  However, it can be useful to know when coarse structures are compatible with other structures that may be present.  The following definition also comes from \cite{HPR}.

\begin{definition} \label{ctt}
Let $X$ be a coarse space.  Then $X$ is called a {\em coarse topological space} if it is equipped with a Hausdorff topology such that every entourage is open, and the closure of every bounded set is compact.
\end{definition}

In fact it is easy to demonstrate that any coarse topological space is locally compact, and the bounded sets are precisely those which are precompact.  If $X$ is a coarse space equipped with some topology, we say that the topology and coarse structure are {\em compatible} when $X$ is a coarse topological space.

The main examples of coarse structures we will use are the bounded coarse structure on a metric space, and the continuously controlled coarse structure arising from a compactification of a topological space.

\begin{definition}
If $X$ is a proper metric space, the {\em bounded coarse structure} is the unital coarse structure formed by defining the entourages to be {\em neighbourhoods of the diagonal}:
$$D_R = \{ (x,y) \in X\times X \ |\ d(x,y)<R \}$$
\end{definition}

\begin{definition}
Let $X$ be a Hausdorff space equipped with a compactification $\overline{X}$.\footnote{This statement means that the space $X$ is included in the compact space $\overline{X}$ as a dense subset.}  Then the {\em continuously controlled} coarse structure is formed by defining the entourages to be open subsets $M\subseteq X\times X$ such that the closure $\overline{M} \subseteq \overline{X}\times \overline{X}$ intersects the set $\partial X\times \partial X$\footnote{The boundary, $\partial X$, is the complement $\overline{X}\backslash X$.} only in the diagonal, and each set of the form $M(x)$ is precompact.
\end{definition}

Both of the above examples are in fact coarse topological spaces.  Making further restrictions, other examples of coarse topological spaces include the `monoidal control spaces' considered in \cite{CM}.

Bounded coarse structures are more general that they might at first appear because of the following result, which is proved in \cite{STY}.

\begin{proposition} \label{boundedgeneral}
Let $X$ be a coarse space.  Suppose that $X$ is {\em countably generated} in the sense that there is a sequence of entourages $(M_n)$ such that every entourage $M$ is contained in a finite composition of the form $M_1 M_2\cdots M_n$.

Then the coarse structure on $X$ is equivalent to the bounded coarse structure arising from some metric.
\noproof
\end{proposition}

A {\em coarse map} is a structure-preserving map between coarse spaces.

\begin{definition}

Let $X$ and $Y$ be coarse spaces.  Then a map $f\co X\rightarrow Y$ is said to be {\em coarse} if:

\begin{itemize}

\item The mapping $f\times f\co X\times X\rightarrow Y\times Y$ takes entourages to subsets of entourages

\item For any bounded subset $B\subseteq Y$ the inverse image $f^{-1}[B]$ is also bounded

\end{itemize}

\end{definition}

The main definitions of coarse geometry are motivated by looking at coarse maps between metric spaces equipped with the bounded coarse structure; see \cite{Roe1}.

We can form the category of all coarse spaces and coarse maps.  We call this category the {\em coarse category}.  We call a coarse map $f\co X\rightarrow Y$ a {\em coarse equivalence} if there is a coarse map $g\co Y\rightarrow X$ such that the composites $g\circ f$ and $f\circ g$ are close to the identities $1_X$ and $1_Y$ respectively.

Let $X$ and $Y$ be coarse spaces, equipped with collections of entourages ${\cal E}(X)$ and ${\cal E}(Y)$ respectively.  Then we define the {\em product} of $X$ and $Y$ to be the Cartesian product $X\times Y$ equipped with the coarse structure defined by forming finite compositions and unions of entourages in the set:
$$\{ M\times N \ |\ M\in {\cal E}(X) , N\in {\cal E}(Y) \}$$

Unfortunately, the above product is not a product in the category-theoretic sense since the projections $\pi_X \co X\times Y\rightarrow X$ and $\pi_Y \co X\times Y\rightarrow Y$ are not in general coarse maps.\footnote{Because the inverse image of a bounded set need not be bounded.}  

\begin{definition}
A {\em generalised ray} is the topological space $[0,\infty )$ equipped with some unital coarse structure compatible with the topology.
\end{definition}

We reserve the notation ${\mathbb R}_+$ for the space $[0,\infty )$ equipped with the bounded coarse structure defined by the metric.

The following definition is a generalisation of the notion of a {\em Lipschitz homotopy} given in \cite{Gr} and used in coarse geometry in \cite{Yu}.

\begin{definition}

Let $X$ and $Y$ be coarse spaces.  Then a {\em coarse homotopy} is a map $F\co X\times R\rightarrow Y$ for some generalised ray $R$ such that:

\begin{itemize}

\item The map $X\times R\rightarrow Y\times R$ defined by writing $(x,t)\mapsto (F(x,t),t)$ is a coarse map.

\item For every bounded set $B\subseteq X$ there is a point $T\in R$ such that the function $(x,t)\mapsto F(x,t)$ is constant if $t\geq T$ and $x\in B$.

\item The function $F(-,\infty)$ defined by the formula:
$$x\mapsto \lim_{t\rightarrow \infty} F(x,t)$$
is a coarse map.

\end{itemize}

\end{definition}

The coarse maps $F(-,\infty )$ and $F(-,0)$ are said to be {\em linked by a coarse homotopy}.  Two coarse functions $f\co X\rightarrow Y$ and $g\co X\rightarrow Y$ are said to be {\em coarsely homotopic} if they are linked by a chain of coarse homotopies.  A coarse map $f\co X\rightarrow Y$ is called a {\em coarse homotopy-equivalence} if there is a coarse map $g\co Y\rightarrow X$ such that the composites $g\circ f$ and $f\circ g$ are coarsely homotopic to the identities $1_X$ and $1_Y$ respectively.

Observe that any two coarsely equivalent maps are coarsely homotopic.  The following result gives us an example of spaces that are coarsely homotopic but not coarsely equivalent.

\begin{proposition} \label{homray}
Let ${\mathbb R}^n$ denote Euclidean space equipped with the bounded coarse structure.  Then the inclusion $i\co {\mathbb R}_+\hookrightarrow {\mathbb R}^n\times {\mathbb R}_+$ defined by the formula $i(s) = (0,s)$ is a coarse homotopy-equivalence.   
\end{proposition}

\proof
Define a coarse map $p\co {\mathbb R}^n\times {\mathbb R}_+\rightarrow {\mathbb R}_+$ by writing $p(x,s) = \| x\| + s$.  The composition $p\circ i$ is equal to the identity $1_{{\mathbb R}_+}$.  We can define a coarse homotopy linking the functions $1_{{\mathbb R}^n\times {\mathbb R}_+}$ and $i\circ p$ by the formula:
$$F(x,s,t) = \left\{ \begin{array}{ll}
(x\cos (\frac{t}{\| x \|} ) , s + \| x \| \sin (\frac{t}{\| x \|} )) & t\leq \frac{\pi \| x \|}{2} \\
(0, s + \| x \| ) & t\geq \frac{\pi \| x\|}{2} \\
\end{array} \right.\eqno{\qed}$$

Actually, the above construction of a homotopy equivalence can be generalised to some situations involving continuous control.  We will come back to this point in section \ref{CW}.

We have seen that there is no good notion of products in the coarse category.  Fortunately, the situation is rather better when we look at coproducts and some more general colimits.

\begin{definition}
Let $\{ X_\lambda \ |\ \lambda \in \Lambda \}$ be a collection of coarse spaces.  Let ${\cal E}_\lambda$ denote the collection of entourages of the space $X_\lambda$, and let ${\cal B}_\lambda$ denote
the collection of bounded subsets.  Then the disjoint union, $\coprod_{\lambda \in \Lambda} X_\lambda$, is equipped with a coarse structure defined by forming the set of entourages:
$$\bigcup_{\lambda \in \Lambda} {\cal E}_\lambda \cup \bigcup_{\lambda , \lambda' \in \Lambda} \{ B\times B' \ |\ B \in {\cal B}_\lambda , B' \in {\cal B}_{\lambda'} \}$$
\end{definition}

It is straightforward to verify that the disjoint union of a collection of coarse spaces defines a coproduct in the coarse category.

\begin{definition}
Let $X$ be a coarse space equipped with an equivalence relation $\sim$ and quotient map $\pi \co X\rightarrow X/\sim$.  Then the quotient $X/\sim$ is equipped with a coarse structure formed by defining the set of entourages to be the collection:
$$\{ \pi [M] \ |\ M\subseteq X\times X \textrm{ is an entourage} \}$$
\end{definition}

In general the quotient map $\pi \co X\rightarrow X/\sim$ is not a coarse map, since the inverse image of a bounded subset need not be bounded.  However, the map $\pi$ is a coarse map when each equivalence class of the relation $\sim$ is finite.

\section{Coarse homology theories} \label{chtsec}

A {\em coarse homology theory} is a functor on the category of coarse spaces that enables us to distinguish different coarse structures.  It is the analogue in the world of coarse geometry of a generalised homology theory in the world of topology.  

\begin{definition} \label{cht}
A {\em coarse homology theory} consists of a collection of functors, $\{ HX_p \}_{p\in {\mathbb Z}}$, from the category of coarse spaces to the category of Abelian groups such that the following axioms hold:

\begin{itemize}

\item Coarse homotopy-invariance:

For any two coarsely homotopic maps $f\co X\rightarrow Y$ and $g\co X\rightarrow Y$, the induced maps $f_\ast \co HX_p (X)\rightarrow HX_p (Y)$ and $g_\ast \co HX_p (X)\rightarrow HX_p (Y)$ are equal.

\item Excision axiom:

Consider a decomposition $X = A\cup B$ of a coarse space $X$.  Suppose that for all entourages $m\subseteq X\times X$ we can find an entourage $M\subseteq X\times X$ such that $m(A)\cap m(B)\subseteq M(A\cap B)$.  Consider the inclusions $i\co A\cap B\hookrightarrow A$, $j\co A\cap B\hookrightarrow B$, $k\co A\hookrightarrow X$, and $l\co B\hookrightarrow X$.  Then we have a natural map $d \co HX_p(X)\rightarrow HX_{p-1}(A\cap B)$ and a long exact sequence:
$$\xymatrix@=8pt{
{} \ar[r] & HX_p (A\cap B) \ar[r]^-{\alpha} & HX_p (A)\oplus HX_p(B) \ar[r]^-{\beta} & HX_p (X) \ar[r]^-{d} & HX_{p-1}(A\cap B) \ar[r] & {}
}$$
where $\alpha = (i_\ast , -j_\ast )$ and $\beta = k_\ast + l_\ast$.

\end{itemize}

The coarse homology theory $\{ HX_p \}$ is said to satisfy the {\em large scale axiom} if the groups $HX_p (+)$ are all trivial, where $+$ denotes the one-point topological space.

\end{definition}

A decomposition, $X = A\cup B$, of a coarse space $X$ is said to be {\em coarsely excisive} if the coarse excision axiom applies, that is to say for all entourages $m\subseteq X\times X$ we can find an entourage $M\subseteq X\times X$ such that $m(A)\cap m(B)\subseteq M(A\cap B)$.  The long exact sequence:
$$\xymatrix@=10pt{
{} \ar[r] & HX_p (A\cap B) \ar[r] & HX_p (A)\oplus HX_p(B) \ar[r] & HX_p (X) \ar[r] & HX_{p-1}(A\cap B) \ar[r] & {}
}$$
is called the {\em coarse Mayer-Vietoris sequence}.

The process of coarsening, described in \cite{HR2} and \cite{Roe1}, is used to construct coarse homology theories on the category of proper metric spaces equipped with their bounded coarse structures.  This process can be significantly generalised.

\begin{definition} \label{good}
Let $X$ be a coarse space.  A {\em good cover} of $X$ is a cover $\{ B_i \ |\ i\in I \}$ such that each set $B_i$ is bounded, and each point $x\in X$ lies in at most finitely many of the sets $B_i$.
\end{definition}

Let $\cal U$ and $\cal V$ be good covers of the space $X$.  Then we say that the cover $\cal V$ is a {\em coarsening} of the cover $\cal U$, and we write ${\cal U}\leq {\cal V}$, if there is a map $\phi \co {\cal U}\rightarrow {\cal V}$ such that $U\subseteq \phi [U]$ for all sets $U\in {\cal U}$.\footnote{A coarsening is the opposite of a refinement, as considered in \v{C}ech cohomology theory.}  The map $\phi$ is called a {\em coarsening map}.

\begin{definition}
A directed family of good covers of $X$, $({\cal U}_i ,\phi_{ij} )_{i\in I}$, is said to be a {\em coarsening family} if there is a family of entourages $(M_i )$ such that:

\begin{itemize}

\item For all sets $U\in {\cal U}_i$ there is a point $x\in X$ such that $U\subseteq M_i (x)$.

\item Let $x\in X$ and suppose that $i<j$.  Then there is a set $U\in {\cal U}_j$ such that $M_i (x)\subseteq U$.

\item Let $M\subseteq X\times X$ be an entourage.  Then $M\subseteq M_i$ for some $i\in I$.

\end{itemize}

\end{definition}

It is proved in \cite{Roe3} that any proper metric space, with its bounded coarse structure, admits a coarsening family.  The following result is a generalisation of this fact.

\begin{proposition} \label{coarsenings}
Let $X$ be a unital coarse topological space.  Then $X$ has a coarsening family.
\end{proposition}

\begin{proof}
Let $\{ M_i \ |\ i\in I \}$ be a generating set for the collection of entourages in the space $X$, in the sense that every entourage $M$ is contained in some finite composite $M_{i_1} M_{i_2} \cdots M_{i_k}$.\footnote{In view of proposition \ref{boundedgeneral} this generating set is probably uncountable.  However, note that we can always find such a generating set simply by considering the set of all entourages.} Without loss of generality, let us assume that each generator $M_i$ contains the diagonal $X\times X$, is {\em symmetric} in the sense that $M_i^{-1} = M_i$, and every finite composite $M_{i_1} M_{i_2} \cdots M_{i_k}$ of generators is an entourage.

Define $J$ to be the set of finite sequences of elements of the set $I$.  The set $J$ can be ordered; we say $j_1 \leq j_2$ if $j_1$ is a subsequence of $j_2$.

Let $i\in I$.  Let us say a subset $S\subseteq X$ is {$M_i$-sparse} if $(x,y)\notin M_i$ for all points $x,y \in S$ such that $x\neq y$.  By Zorn's lemma there is a maximal $M_i$-sparse set $A_i$.

Consider a finite sequence $j = ( i_1 ,\ldots , i_k ) \in J$.   
Form the composite $M_j = M_{i_1}\ldots M_{i_k}$ and union $A_j = A_{i_1}
\cup \cdots \cup A_{i_k}$, and look at the collection of bounded sets:
\[ {\cal U}_j = \{ M_j (x) \ |\ x\in A_j \} \]
%
Since each of the generators is symmetric and contains the diagonal, the composite $M_j$ contains each of the generators that form it.  By maximality of the $M_{i_1}$-sparse set $A_{i_1}$, the collection ${\cal U}_j$ is a cover of the space $X$.  Choose a point $x_0 \in X$.  Then the set
\[ \{ x\in A_j \ |\ x_0 \in M_j (x) \} \]
is a discrete subset of a compact set, and is therefore finite.  Thus the cover ${\cal U}_j$ is good in the sense of definition \ref{good}.

We now need to check that the family of covers $( {\cal U}_j )$ forms a coarsening family.  Consider the family of entourages $(M_j)$.  By construction, every set $U\in {\cal U}_j$ takes the form $M_j (x)$ for some point $x\in X$, and every entourage $M$ is contained in some entourage of the form $M_j$.

All that remains is to check that for every point $x\in X$ and every pair of elements $j,j' \in J$ such that $j<j'$ there is a set $U\in {\cal U}_{j'}$ such that $M_j (x)\subseteq U$.

It suffices to look at the special case where we have sequences $j = (i_1 , \ldots , i_k)$ and $j' = (i_1 , \ldots , i_k ,i_{k+1})$ in the indexing set $J$.  By maximality of the $M_{i_{k+1}}$-sparse set $A_{i_{k+1}}$ we can find a point $y\in A_{i_{k+1}}$ such that $x\in M_{i_{k+1}}(y)$.  Hence $M_j (x)\subseteq M_{j'}(y)$.  
\end{proof}

In particular, note that any generalised ray admits a coarsening family.

Recall that the {\em nerve}, $|{\cal U}_i |$, of a cover ${\cal U}_i$ is a simplicial set with one vertex for every set in the cover.  Vertices represented by sets $U_1 , \ldots , U_n\in {\cal U}_i$ are spanned by an $n$-simplex if and only if the intersection $U_1 \cap \cdots \cap U_n$ is non-empty.  If $( {\cal U}_i , \phi_{ij} )$ is a coarsening family, then each map $\phi_{ij} \co {\cal U}_i\rightarrow {\cal U}_j$ induces a map of simplicial sets ${\phi_{ij}}_\ast \co |{\cal U}_i |\rightarrow |{\cal U}_j |$.  The map ${\phi_{ij}}_\ast$ is proper because the covers we consider are locally finite.\footnote{That is to say each point of our space lies in only finitely many sets of a particular cover.}

\begin{definition}
Let $\{ H^\mathrm{lf}_p \}$ be a generalised locally finite homology theory on the category of simplicial sets.  Let $({\cal U}_i , \phi_{ij})$ be a coarsening family on the coarse space $X$.  Then we write:
$$HX_p (X) = \lim_{\rightarrow \atop i} H_p^\mathrm{lf} |{\cal U}_i |$$
\end{definition}

As the above terminology suggests, it will turn out that the assignment $X\mapsto HX_p (X)$ forms a coarse homology theory.

\begin{proposition} \label{func}
Let $X$ and $Y$ be coarse spaces equipped with coarsening families $({\cal U}_i )_{i\in I}$ and $({\cal V}_j)_{j\in J}$ respectively.  Let $f\co X\rightarrow Y$ be a coarse map.  Then there is a functorially induced homomorphism $f_\ast \co HX_p (X)\rightarrow HX_p (Y)$
\end{proposition}

\begin{proof}
Since the map $f$ is a coarse map, it takes entourages to entourages.  The given conditions on coarsening families ensure that for each element $i\in I$, $f[{\cal U}_i ]\leq {\cal V}_j$ for some element $j\in J$.

Hence we have an induced map of simplicial sets $f_\ast \co |{\cal U}_i |\rightarrow |{\cal V}_j|$.  This map is proper since the map $f$ is coarse.  By definition of the nerve of a cover, up to proper homotopy the map $f_\ast$ is independent of any choice of coarsening family.  We therefore obtain a homomorphism:
$$f_\ast \co \lim_{\rightarrow \atop i} H_p^\mathrm{lf}|{\cal U}_i |\rightarrow \lim_{\rightarrow \atop j} H_p^\mathrm{lf}|{\cal V}_j |$$
by taking locally finite homology groups followed by the direct limit.
\end{proof}

\begin{corollary}
The group $HX_p (X)$ does not, up to isomorphism, depend on the choice of coarsening family.
\noproof
\end{corollary}

\begin{lemma} \label{hominvar}
Let $f\co X\rightarrow Y$ and $g\co X\rightarrow Y$ be coarsely homotopic maps between coarse spaces that admit coarsening families.  Then the induced maps $f_\ast \co HX_p (X)\rightarrow HX_p (Y)$ and $g_\ast \co HX_p (X)\rightarrow HX_p (Y)$ are equal.
\end{lemma}

\begin{proof}
Without loss of generality suppose we have a coarse homotopy $F\co X\times R\rightarrow Y$ such that $f = F(-,0)$ and $g = F(- , \infty )$.  Let $({\cal U}_i)_{i\in I}$, $({\cal V}_j)_{j\in J}$, and $({\cal W}_k )_{k\in K}$ be coarsening families for the spaces $X$, $R$, and $Y$ respectively.  Then for all elements $i\in I$ and $j\in J$ we can find an element $k\in K$ such that $F({\cal U}_i \times {\cal V}_j )\leq {\cal W}_k$.  Hence we have an induced proper map of simplicial sets:
$$F_\ast \times 1 \co |{\cal U}_i |\times |{\cal V}_j |\rightarrow |{\cal W}_k |\times |{\cal V}_j |$$

Now, let us write the cover ${\cal V}_j$ of the generalised ray $R$ as a sequence of bounded sets $(V_n)$ where $\sup V_{n+1}\geq \sup V_n$ and $V_{n+1}\cap V_n \neq \emptyset$.  Then we can define a continuous map $[0,\infty )\rightarrow |{\cal V}_j |$ by sending a natural number $n\in {\mathbb N}$ to the vertex $V_n$ and a point $t\in (n , n+1)$ to the appropriate point on the edge joining the vertices $V_n$ and $V_{n+1}$.  Hence we have a proper continuous map
$$F_\ast \times 1\co |{\cal U}_i |\times [0,\infty )\rightarrow |{\cal W}_k |\times [0,\infty )$$
such that for each point $x\in |{\cal U}_i|$ the function $t\mapsto F_\ast (x,t)$ is eventually constant.  

The induced maps $f_\ast \co |{\cal U}_i |\rightarrow |{\cal W}_k |$ and $g_\ast \co |{\cal U}_i |\rightarrow |{\cal W}_k |$ are thus properly homotopic; we obtain the appropriate proper homotopy by renormalising the map $F_\ast \co |{\cal U}_i |\times [0,\infty )\rightarrow |{\cal W}_k |$.  The maps $f_\ast$ and $g_\ast$ therefore induce the same map at the level of locally finite homology.  Taking direct limits, the maps $f_\ast \co HX_p (X)\rightarrow HX_p (Y)$ and $g_\ast \co HX_p (X)\rightarrow HX_p (Y)$ must be equal.
\end{proof}

\begin{lemma} \label{excision}
Let $X$ be a coarse space that admits some coarsening family.  Suppose we have a coarsely excisive decomposition $X = A\cup B$ and inclusions $i\co A\cap B\hookrightarrow A$, $j\co A\cap B\hookrightarrow B$, $k\co A\hookrightarrow X$, and $l\co B\hookrightarrow X$.  Then we have a natural map $d \co HX_p(X)\rightarrow HX_{p-1}(A\cap B)$ and a long exact sequence:
$$\xymatrix@=10pt{
{} \ar[r] & HX_p (A\cap B) \ar[r]^-{\alpha} & HX_p (A)\oplus HX_p(B) \ar[r]^-{\beta} & HX_p (X) \ar[r]^-{d} & HX_{p-1}(A\cap B) \ar[r] & {}
}$$
where $\alpha = (i_\ast , -j_\ast )$ and $\beta = k_\ast + l_\ast$.
\end{lemma}

\begin{proof}
Let $({\cal U}_i )_{i\in I}$ be a coarsening family for the space $X$.  Then the spaces $A$ and $B$ have coarsening families $({\cal U}_i |_A)$ and $({\cal U}_i |_B)$ defined by writing:
$${\cal U}_i |_A = \{ U\cap A \ |\ U\in {\cal U}_i , U\cap A\neq \emptyset \} \qquad {\cal U}_i |_B = \{ U\cap B \ |\ U\in {\cal U}_i , U\cap B\neq \emptyset \}$$ 
respectively.

Write $A_i = | {\cal U}_i |_A |$, $B_i = | {\cal U}_i |_B |$, and $X_i = |{\cal U}_i |$.  The decomposition $X = A\cup B$ is coarsely excisive, so if we look at interiors then:
$$A_i^0 \cup B_i^0 = X_i$$

By the existence of Mayer-Vietoris sequences in ordinary homology (see for example \cite{Sp}) we have natural maps $d \co H^\mathrm{lf}_p( X_i )\rightarrow H^\mathrm{lf}_{p-1}(A_i \cap B_i)$ and exact sequences:
$$\xymatrix@=10pt{
{} \ar[r] & H^\mathrm{lf}_p (A_i\cap B_i) \ar[r] & H^\mathrm{lf}_p (A_i)\oplus H^\mathrm{lf}_p(B_i) \ar[r] & H^\mathrm{lf}_p (X_i) \ar[r]^-{d} & H^\mathrm{lf}_{p-1}(A_i\cap B_i) \ar[r] & {}
}$$

Taking direct limits, we have an exact sequence:
$$\xymatrix@=10pt{
{} \ar[r] & HX_p (A\cap B) \ar[r]^-{\alpha} & HX_p (A)\oplus HX_p(B) \ar[r]^-{\beta} & HX_p (X) \ar[r]^-{d} & HX_{p-1}(A\cap B) \ar[r] & {}
}$$
as required.
\end{proof}

Lemma \ref{hominvar} and lemma \ref{excision} together imply the following result.

\begin{theorem}
Let $\{ H^\mathrm{lf}_p \}$ be a locally finite homology theory.  Then the collection of functors $\{ HX_p \}$ defines a coarse homology theory on the category of coarse topological spaces.
\noproof
\end{theorem}

The coarse homology theories considered above {\em do not} satisfy the large scale axiom.

Now, let $X$ be a coarse paracompact topological space equipped with an open coarsening family $({\cal U}_i )$.  Let $\{ \varphi_U \ |\ U\in {\cal U}_i \}$ be a partition of unity subordinate the the open cover ${\cal U}_i$.  Given a point $x\in X$ there are only finitely many sets $U\in {\cal U}_i$ such that $x\in U$.  The sum $\sum_{U\in {\cal U}_i } \varphi_U (x) U$ represents a point in the simplex spanning the vertices represented by these sets.  We can therefore define a proper continuous map $\kappa \co X\rightarrow |{\cal U}_i |$ by the formula:
$$\kappa_i (x) = \sum_{U\in {\cal U}_i } \varphi_U (x) U$$

Considering locally finite homology and taking direct limits we obtain a {\em coarsening map}:
$$c\co H_n^\mathrm{lf}(X)\rightarrow HX_n (X)$$

\begin{proposition}
The coarsening map $c$ does not depend upon the choice of partition of unity.
\end{proposition}

\begin{proof}
Let $\{ \varphi_U \ |\ U\in {\cal U}_i \}$ and $\{ \varphi_U' \ |\ U\in {\cal U}_i \}$ be partitions of unity subordinate to the open cover ${\cal U}_i$.  Then the resulting maps $\kappa_i$ and $\kappa'_i$ are properly homotopic, and so induce the same maps at the level of locally finite homology.
\end{proof} 

It is proved in \cite{HR2} that the coarsening map, $c$, is an isomorphism whenever $X$ is a proper metric space\footnote{Equipped with the bounded coarse structure.} with good local properties.  In particular, it is easy to see that the map $c$ is an isomorphism if the space $X$ is a single point.

\section{Relative coarse homology}

Homology theories are usually defined by looking at pairs of topological spaces.  There is a corresponding notion in the coarse setting.

\begin{definition}
A {\em pair}, $(X,A)$, of coarse spaces consists of a coarse space $X$ and a subspace $A\subseteq X$.  A {\em coarse map of pairs} $f\co (X,A)\rightarrow (X,B)$ is a coarse map $f\co X\rightarrow Y$ such that $f[A]\subseteq B$.  A {\em relative coarse homotopy} $F\co (X,A)\times R\rightarrow (Y,B)$ is a coarse homotopy $F\co X\times R\rightarrow Y$ such that $F(a,t)\in B$ for all points $a\in A$ and $t\in R$.  Two coarse maps of pairs $f,g\co (X,A)\rightarrow (Y,B)$ are {\em relatively coarsely homotopic} if they are linked by a chain of relative coarse homotopies.
\end{definition}

We now present the coarse version of the Eilenberg-Steenrod axioms.

\begin{definition}
A {\em relative coarse homology theory} consists of a collection of functors, $\{ HX_p \}_{p\in {\mathbb Z}}$, from the category of pairs of coarse spaces to the category of Abelian groups, together with natural transformations $\partial \co HX_p (X,A)\rightarrow HX_p (A,\emptyset )$ such that the following axioms hold:

\begin{itemize}

\item Coarse homotopy-invariance:

Let $f\co (X,A)\rightarrow (Y,B)$ and $g\co (X,A)\rightarrow (Y,B)$ be relatively coarsely homotopic coarse maps.  Then the induced maps $f_\ast \co HX_p (X,A)\rightarrow HX_p (Y,B)$ and $g_\ast \co HX_p (X,A)\rightarrow HX_p (Y,B)$ are equal.

\item Long exact sequence axiom:

The inclusions $i\co (A,\emptyset )\hookrightarrow (X,\emptyset )$ and $j\co (X,\emptyset )\hookrightarrow (X,A)$ induce a long exact sequence:
$$\longrightarrow HX_p (A,\emptyset) \stackrel{i_\ast}{\longrightarrow} HX_p (X,\emptyset ) \stackrel{j_\ast}{\longrightarrow} HX_p (X,A)\stackrel{\partial}{\longrightarrow} HX_{p-1}(A,\emptyset )\longrightarrow$$

\item Excision axiom:

Suppose we have a coarsely excisive decomposition $X = A\cup B$.  Then the inclusion $(A,A\cap B)\hookrightarrow (X,B)$ induces an isomorphism $HX_p (A,A\cap B)\rightarrow HX_p (X,B)$.

\item Large scale axiom:

Let $+$ denote the one-point coarse space.  Then the groups $HX_p (+ ,\emptyset )$ are all trivial.

\end{itemize}

\end{definition}

Let $\{ HX_p \}$ be a relative coarse homology theory.  Then we write $HX_p (X) = HX_p (X, \emptyset )$.  The assignment $X\mapsto HX_p (X)$ is a coarse homology theory in the sense of definition \ref{cht} because of the following result.

\begin{proposition}
Let $X$ be a coarse space, equipped with a coarsely excisive decomposition $X = A\cup B$ and inclusions $i\co A\cap B\hookrightarrow A$, $j\co A\cap B\hookrightarrow B$, $k\co A\hookrightarrow X$, and $l\co B\hookrightarrow X$.  Then we have a natural map $d \co HX_p(X)\rightarrow HX_{p-1}(A\cap B)$ and a long exact sequence:
$$\xymatrix@=10pt{
{} \ar[r] & HX_p (A\cap B) \ar[r]^-{\alpha} & HX_p (A)\oplus HX_p(B) \ar[r]^-{\beta} & HX_p (X) \ar[r]^-{d} & HX_{p-1}(A\cap B) \ar[r] & {}
}$$
where $\alpha = (i_\ast , -j_\ast )$ and $\beta = k_\ast + l_\ast$.
\end{proposition}

\begin{proof}
We have a commutative diagram of exact sequences:
$$\xymatrix{
HX_p (A\cap B) \ar[r]^-{i_\ast} \ar[d]_{j_\ast} & HX_p (A) \ar[r] \ar[d]_{k_\ast} & HX_p (A,A\cap B) \ar[r]^-{\partial} \ar[d]^-{c} & HX_{p-1}(A\cap B) \ar[d]^-{j_\ast} \\
HX_p (B) \ar[r]^-{l_\ast} & HX_p (X) \ar[r]^-{b} & HX_p (X,B) \ar[r]^-{\partial} & HX_{p-1}(B) \\
}$$

The excision map $c\co HX_p (A,A\cap B)\rightarrow HX_p (X,B)$ is an isomorphism.  Write $d = \partial c^{-1}b$.  Then a diagram chase tells us that we have a long exact sequence:
$$\xymatrix@=10pt{
{} \ar[r] & HX_p (A\cap B) \ar[r]^-{\alpha} & HX_p (A)\oplus HX_p(B) \ar[r]^-{\beta} & HX_p (X) \ar[r]^-{d} & HX_{p-1}(A\cap B) \ar[r] & {}
}$$
as desired.
\end{proof}

Now let $X$ be a subspace of the sphere $S^{n-1}$.  The {\em open cone} on $X$, ${\cal O}X$, is the subset of Euclidean space ${\mathbb R}^n$ consisting of points:
$$\{ tx \ |\ t\in (0 , \infty ) , x\in X \}$$

The coarse geometry of the cone ${\cal O}X$ is closely related to the topology of the space $X$--- see for example \cite{HR2}, \cite{PRW}, and \cite{Roe1} for instances of this phenomenon.

\begin{definition}
Let $X$ be a subspace of the sphere $S^{n-1}$.   Let $r\co (0,\infty )\rightarrow (0, \infty )$ be any map such that $r(t)\rightarrow \infty$ as $t\rightarrow \infty$ and the inequality:
$$|r(s) - r(t)|\leq |s-t|$$
holds for all points $s,t\in {\mathbb R}^{\geq 0}$.  Then the {\em radial contraction} associated to $r$ is the map $\rho \co {\cal O}X\rightarrow {\cal O}X$ defined by the formula:
$$\rho (tx) = r(t) x$$
\end{definition}

The following result is proved in \cite{HR2}.

\begin{lemma} \label{homcoarse}
Let $X$ be a subspace of the sphere $S^{n-1}$ and let $Y$ be any metric space, equipped with the bounded coarse structure arising from the metric.  Then for any continuous map $f\co {\cal O}X\rightarrow Y$ there is a radial contraction $\rho \co {\cal O}X\rightarrow {\cal O}X$ such that the composite $f \circ \rho \co {\cal O}X\rightarrow Y$ maps entourages to entourages $($with respect to the bounded coarse structures$)$.
\noproof
\end{lemma}

It is easy to see that any radial contraction is a coarse map, and is coarsely homotopic to the identity map.

\begin{theorem}
Let $HX_p$ be a relative coarse homology theory.  Then we can define a generalised homology theory on the category of pairs of subspaces of spheres by writing:
$$H_p (X,A) = HX_p ({\cal O}X ,{\cal O}A)$$
whenever $(X,A)$ is a pair of subspaces of some sphere $S^n$.
\end{theorem}

\begin{proof}
If $f\co (X,A)\rightarrow (Y,B)$ is a continuous map of pairs there is an induced continuous map ${\cal O}f \co ({\cal O}X,{\cal O}A)\rightarrow ({\cal O}Y,{\cal O}B)$.  By lemma \ref{homcoarse} we can define a coarse map ${\cal O}f\circ \rho \co ({\cal O}X , {\cal O}A)\rightarrow ({\cal O}Y , {\cal O}B)$ for some radial contraction $\rho$.  The coarse homotopy type of the composition $f\circ \rho$ does not depend on the choice of contraction $\rho$ so we obtain a functorially induced map $f_\ast \co H_p (X,A)\rightarrow H_p (Y,B)$.

Let $F\co (X,A)\times [0,1]\rightarrow (Y,B)$ be a relative homotopy.  Choose a radial contraction $\rho \co {\cal O}X\rightarrow {\cal O}X$ such that each map ${\cal O}F(-,t)\circ \rho \co ({\cal O}X ,{\cal O}A)\rightarrow ({\cal O}Y ,{\cal O}B)$ is a coarse map.  Then we can define a coarse homotopy $G\co ({\cal O}X , {\cal O}A)\times {\mathbb R}_+\rightarrow ({\cal O}Y, {\cal O}B)$ by the formula:
$$G(tx, s) = \left\{ \begin{array}{ll}
F(\rho (tx) , st^{-1} ) & 0\leq s<t \\
F(\rho (tx) , 1) & s\geq t \\
\end{array} \right.$$
for all points $x\in X$, $s\in {\mathbb R}_+$, and $t\in (0,\infty )$.  Hence the maps $\rho \circ {\cal O}F(-,0)$ and $\rho \circ {\cal O}F(-,1)$ are relatively coarsely homotopic, so the induced maps $F(-,0)_\ast$ and $F(-,1)_\ast$ are equal.

If $(X,A)$ is a pair of subspaces of the sphere $S^{n-1}$ then $({\cal O}X,{\cal O}A)$ is a pair of coarse spaces so we have natural maps $\partial \co H_p (X,A)\rightarrow H_{p-1}(A,\emptyset )$ and a long exact sequence 
$$\xymatrix@1{
{} \ar[r] & H_p (A,\emptyset) \ar[r]^-{i_\ast} & H_p (X,\emptyset ) \ar[r]^-{j_\ast} & H_p (X,A) \ar[r]^-{\partial} & H_{p-1}(A,\emptyset ) \ar[r] & {}
}$$
where the maps $i_\ast$ and $j_\ast$ are induced by the inclusions $i\co (A,\emptyset )\hookrightarrow (X,\emptyset )$ and $j\co (X,\emptyset )\hookrightarrow (X,A)$ respectively.

Finally, suppose we have a pair $(X,A)$ and a subset $K\subseteq A$ such that $\overline{K}\subseteq A^0$.  Since the space $\overline{K}$ is compact we can find a real number $\varepsilon >0$ such that the neighbourhood $N_\varepsilon (K)$ is a subset of the space $A$.  Hence for all $R>0$ there exists $S>0$ such that:
$$N_R(A)\cap N_R(X\backslash K )\subseteq N_S(A\cap (X\backslash K))$$

Looking at cones, we see that the decomposition ${\cal O}X = {\cal O}A \cup {\cal O}(X\backslash K)$ is coarsely excisive, and the excision axiom follows.  This completes the proof.
\end{proof}

\section{Coarse $CW$-complexes}  \label{CW} 

The cone of the sphere $S^{n-1}$ is the Euclidean space ${\mathbb R}^n$ and the cone of the ball $D^n$ is the half-space ${\mathbb R}^n\times {\mathbb R}^{\geq 0}$.  By proposition \ref{homray} the coarse space ${\mathbb R}^n\times {\mathbb R}_+$ is coarsely homotopy-equivalent to the ray ${\mathbb R}_+$.

This prompts the following definition.

\begin{definition}
A {\em coarse $0$-cell} is a generalised ray.  A {\em coarse $n$-cell} is the half-space ${\mathbb R}^n\times [0,\infty )$ equipped with some unital coarse structure compatible with the topology such that the inclusion $i\co [0,\infty )\hookrightarrow {\mathbb R}^n \times [0,\infty )$ defined by the formula $i(s) = (0,s)$ is a coarse homotopy-equivalence.
\end{definition}

If $DX^n$ is the space ${\mathbb R}^n\times [0,\infty )$ equipped with some coarse structure that makes it a coarse $n$-cell, then we refer to the space $SX^{n-1} = \{ (x,0) \ |\ x\in {\mathbb R}^n \}$ as the {\em boundary} of the coarse cell $DX^n$.  

\begin{definition}
A {\em coarse $(n-1)$-sphere} is a boundary of some coarse $n$-cell.
\end{definition}

We have already seen that the space ${\mathbb R}^n \times {\mathbb R}_+$, equipped with the bounded coarse structure, is a coarse $n$-cell.  However, other coarse $n$-cells are possible if we look at continuous control.

\begin{proposition}
Let $R$ be a generalised ray, with the coarse structure defined by looking at continuous control at infinity with respect to the one point compactification of the space $[0,\infty )$.  Then the coarse space $(R\coprod R)^n \times R$ is a coarse $n$-cell.
\end{proposition}

\begin{proof}
The topological space ${\mathbb R}^n \times [0,\infty )$ can be compactified by adding a `hemisphere at infinity'.  The coarse space $(R\coprod R)^n \times R$ can be viewed as the space ${\mathbb R}^n \times [0,\infty )$ with coarse structure defined by looking at continuous control with respect to this compactification.

As in proposition \ref{homray} we can define a coarse map $p\co {\mathbb R}^n\times {\mathbb R}_+\rightarrow {\mathbb R}_+$ by writing $p(x,s) = \| x\| + s$.  The composition $p\circ i$ is equal to the identity $1_{{\mathbb R}_+}$, and we can define a coarse homotopy linking the functions $1_{{\mathbb R}^n\times {\mathbb R}_+}$ and $i\circ p$ by the formula:
$$F(x,s,t) = \left\{ \begin{array}{ll}
(x\cos (\frac{t}{\| x \|} ) , s + \| x \| \sin (\frac{t}{\| x \|} )) & t\leq \frac{\pi \| x \|}{2} \\
(0, s + \| x \| ) & t\geq \frac{\pi \| x\|}{2} \\
\end{array} \right.$$
\end{proof}

\begin{definition}
Let $(X,A)$ be a pair of coarse spaces, and let $f\co A\rightarrow Y$ be a continuous map.  Then we define the space $X\cup_f Y$ to be the quotient:
$$X\cup_f Y = \frac{X\coprod Y}{a\sim f(a)}$$
\end{definition}

\begin{proposition} \label{deform}
Let $\pi_X \co X\rightarrow X\cup_f Y$ and $\pi_Y \co Y\rightarrow X\cup_f Y$ be the canonical maps associated to the quotient $X\cup_f Y$.  Then the images $\pi_X [X]$ and $\pi_Y [Y]$ are coarsely equivalent to the original spaces $X$ and $Y$ respectively.  Further, we have a coarsely excisive decomposition $X\cup_f Y = \pi_X [X]\cup \pi_Y [Y]$.
\end{proposition}

\begin{proof}
The maps $\pi_X$ and $\pi_Y$ are injective, and the coarse structures of the images $\pi_X [X]$ and $\pi_Y [Y]$ are defined to be those inherited from the spaces $X$ and $Y$, respectively, under these maps.  This establishes the first part of the proposition.

Certainly the space $X\cup_f Y$ is equal to the union $\pi_X [X]\cup \pi_Y [Y]$.  Let $m$ be an entourage for the space $X\cup_f Y$.  We want to find an entourage $M$ such that:
$$m(\pi_X [X])\cap m(\pi_Y [Y]) \subseteq M(\pi_X [X]\cap \pi_Y [Y])$$

By definition of the coarse structure on the space $X\cup_f Y$ we can write:
$$m = \pi [ m_X \cup m_Y \cup B_X\times B_Y' \cup B_Y\times B_X']$$
where $m_X$ and $m_Y$ are entourages for the spaces $X$ and $Y$ respectively, $B_X , B_X'\subseteq X$ and $B_Y ,B_Y'\subseteq Y$ are bounded subsets, and $\pi \co X\coprod Y\rightarrow X\cup_f Y$ is the canonical quotient map.  Observe that:
$$m(\pi_X [X])\cap m(\pi_Y [Y]) \subseteq (\pi_X [X]\cap \pi_Y [Y])\cup \pi_X [B_X] \cup \pi_Y [B_Y]$$

Let $D_X$ and $D_Y$ be the diagonals in the spaces $X$ and $Y$ respectively.  Consider any point $p\in A$ and let $M$ be an entourage containing the image:
$$\pi [D_X \cup D_Y \cup (B_X\times p) \cup (B_Y\times p) ]$$

Then:
$$m(\pi_X [X])\cap m(\pi_Y [Y]) \subseteq M(\pi_X [X]\cap \pi_Y [Y])$$
and we are done.
\end{proof}

Let $DX^n$ be a coarse $n$-cell, with boundary $SX^{n-1}$.  Consider a coarse map $f\co SX^{n-1}\rightarrow Y$.  Then the coarse space $DX^n \cup_f Y$ is called the coarse space obtained from $Y$ by {\em attaching a coarse $n$-cell} through the map $f$.

\begin{definition}
A {\em finite coarse $CW$-complex} is a coarse space $X$ together with a sequence
$$X^0 \subseteq X^1 \subseteq \cdots \subseteq X^n = X$$
of subspaces such that:

\begin{itemize}

\item The space $X^0$ is a finite disjoint union of generalised rays.

\item The space $X^k$ is coarsely equivalent to a space obtained by attaching a finite number of coarse $k$-cells to the space $X^{k-1}$.

\end{itemize}

\end{definition}

Assuming that $X^n\backslash X^{n-1} \neq \emptyset$, the number $n$ is called the {\em dimension} of the finite coarse $CW$-complex $X$.

The main purpose of this section is to prove that the axioms determine a coarse homology theory completely on the category of spaces coarsely homotopy-equivalent to finite coarse $CW$-complexes once we know what the coarse homology of a generalised ray and a one-point set is.  

\begin{lemma} \label{deformsequence}
Let $DX^n$ be a coarse $n$-cell with boundary $SX^{n-1}$.  Suppose that the cell $DX^n$ is coarsely homotopy-equivalent to a generalised ray $R$.  Form the space $X = DX^n \cup_f Y$ for some coarse map $f\co SX^n\rightarrow Y$.  Then we have a Mayer-Vietoris sequence
$$\xymatrix@=10pt{
{} \ar[r] & HX_p (SX^{n-1}) \ar[r] & HX_p (Y)\oplus HX_p (R) \ar[r] & HX_p (X) \ar[r] & HX_{p-1}(SX^{n-1}) \ar[r] & {}
}$$
\end{lemma}

\begin{proof}
Let $R'$ denote the space $[1,\infty )$ equipped with the coarse structure inherited from the generalised ray $R$ and let $C$ be the space $\{ (x,t)\in DX^n \ |\ t\geq 1 \}$.  Then the space $C$ is coarsely equivalent to the space $DX^n$, and therefore coarsely homotopy-equivalent to the generalised ray $R$.

Let $B = (SX^{n-1}\times [0,2])\cup_f Y$.  Then the space $B$ is coarsely equivalent to the space $Y$.  The space $X$ can be written as the union $X = B\cup C$, and the intersection $B\cap C$ is coarsely equivalent to the coarse sphere $SX^{n-1}$.  Hence by proposition \ref{deform} we have a Mayer-Vietoris sequence
$$\xymatrix@=10pt{
{} \ar[r] & HX_p (SX^{n-1}) \ar[r] & HX_p (Y)\oplus HX_p (R) \ar[r] & HX_p (X) \ar[r] & HX_{p-1}(SX^{n-1}) \ar[r] & {}
}$$
as required.
\end{proof}

The proof of our main result is now virtually identical to the proof of the corresponding result in classical algebraic topology.\footnote{See, for example, section 8, chapter 4, of \cite{Sp} for a proof of the corresponding classical result.}

\begin{definition}
A {\em natural transformation} between coarse homology theories $\{ HX_p \}$ and $\{ HX'_p \}$ is a sequence of natural transformations $\tau \co HX_p \rightarrow HX_p'$ that takes Mayer-Vietoris sequences appearing in the coarse homology theory $\{ HX_p \}$ to Mayer-Vietoris sequences appearing in the coarse homology theory $\{ HX_p' \}$.
\end{definition}

\begin{lemma} \label{spheres}
Let $\tau \co HX_p\rightarrow HX_p'$ be a natural transformation of coarse homology theories such that the map $\tau \co HX_p (X)\rightarrow HX_p'(X)$ is an isomorphism whenever the space $X$ is a single point or a generalised ray.  Then the map $\tau \co HX_p (X)\rightarrow HX_p'(X)$ is an isomorphism whenever the space $X$ is a coarse sphere.
\end{lemma}

\begin{proof}
Let us write a given coarse $0$-sphere, $SX^0$, as a coarsely excisive union, $SX^0 = R_1\cup R_2$, of two generalised rays.  The intersection, $R_1\cap R_2$, is bounded and is therefore coarsely equivalent to a single point, $+$.  Considering Mayer-Vietoris sequences we have a commutative diagram:
$$\xymatrix@=10pt{
{} \ar[r] & {HX_p (+)}\ar[r] \ar[d] & {HX_p(R_1)\oplus HX_p (R_2)}\ar[r] \ar[d] & {HX_p (SX^0)}\ar[r] \ar[d] & {HX_{p-1}(+)} \ar[r] \ar[d] & {} \\
{} \ar[r] & {HX_p' (+)}\ar[r] & {HX_p'(R_1)\oplus HX_p'(R_2)}\ar[r] & {HX_p' (SX^0)}\ar[r] & {HX_{p-1}'(+)} \ar[r] & {} \\
}$$

The rows in the above diagram are exact.  With the possible exception of the map $\tau \co HX_p (SX^0)\rightarrow HX_p (SX^0)$, the vertical arrows are isomorphisms.  Hence the map $\tau \co HX_p (SX^0)\rightarrow HX_p'(SX^0)$ is also an isomorphism by the five lemma.

Now suppose that the map $\tau$ is an isomorphism for every coarse $(n-1)$-sphere.  Let $S$ be a coarse $n$-sphere.  Then we can write $S$ as a coarsely excisive union, $D_1\cup D_2$, of two $n$-cells, with intersection coarsely equivalent to some coarse $(n-1)$-sphere, $S_0$.  The result now follows by induction if we look at Mayer-Vietoris sequences and apply the five lemma as above.
\end{proof}

\begin{theorem} \label{uniqueness}
Let $\tau \co \{ HX_p \} \rightarrow \{ HX_p' \}$ be a natural transformation of coarse homology theories such that the map $\tau \co HX_p (X)\rightarrow HX_p' (X)$ is an isomorphism whenever the space $X$ is a single point or a generalised ray.  Then the map $\tau \co HX_p (X)\rightarrow HX_p' (X)$ is an isomorphism whenever the space $X$ is coarsely homotopy-equivalent to a finite coarse $CW$-complex.
\end{theorem}

\begin{proof}
The map $\tau \co HX_p (X)\rightarrow HX_p'(X)$ is certainly an isomorphism if the space $X$ is a coarse $CW$-complex with just one cell.  Suppose that the map is an isomorphism whenever the space $X$ is a coarse $CW$-complex with fewer than $m$ cells.

Let $X$ be a coarse $CW$-complex of dimension $n$ that has $m$ cells.  Let $C$ be a cell of dimension $n$ and let $B$ be a coarse $CW$-complex such that the $CW$-complex $X$ is obtained from $B$ by attaching the cell $C$.  By lemma \ref{deformsequence} we have a commutative diagram:

$$\xymatrix@=10pt{
{} \ar[r] & {HX_p(S)} \ar[r] \ar[d] & {HX_p (B)\oplus HX_p (R)} \ar[r] \ar[d] & {HX_p (X)} \ar[r] \ar[d] & {HX_{p-1} (S)} \ar[r] \ar[d] & {} \\
{} \ar[r] & {HX'_p(S)} \ar[r] & {HX'_p (B)\oplus HX'_p (R)} \ar[r] & {HX'_p (X)} \ar[r] & {HX'_{p-1}(S)} \ar[r] & {} \\
}$$
where $S$ is a coarse $(n-1)$-sphere, $R$ is a generalised ray coarsely homotopy-equivalent to the cell $C$, and the rows are Mayer-Vietoris sequences.  All of the vertical arrows except for possibly the map $\tau \co HX_p (X)\rightarrow HX_p'(X)$ are isomorphisms by inductive hypothesis and lemma \ref{spheres}.  Hence by the five lemma the map $\tau \co HX_p(X)\rightarrow HX_p'(X)$ is also an isomorphism.

Thus the result holds for any finite coarse $CW$-complex by induction.  By coarse homotopy-invariance the result follows for any space coarsely homotopy-equivalent to a finite $CW$-complex.
\end{proof}

Now, let $\{ HX_p \}$ and $\{ HX_p' \}$ be {\em relative} coarse homology theories.  Then a {\em natural transformation} between these theories is a sequence of natural transformations $\tau \co HX_p \rightarrow HX_p'$ that takes the long exact sequences appearing in the homology theory $\{ HX_p \}$ to the long exact sequences appearing in the homology theory $\{ HX_p' \}$.  A natural transformation of relative coarse homology theories induces a natural transformation of the corresponding coarse homology theories.

\begin{corollary}
Let $\tau \co \{ HX_p \}\rightarrow \{ HX_p' \}$ be a natural transformation of relative coarse homology theories such that the map $\tau \co HX_p (R)\rightarrow HX_p' (R)$ is an isomorphism for every generalised ray, $R$.  Then the map $\tau \co HX_p (X,A)\rightarrow HX_p' (X,A)$ is an isomorphism for every pair of coarse spaces, $(X,A)$, that are coarsely homotopy-equivalent to finite coarse $CW$-complexes.
\end{corollary}

\begin{proof}
By theorem \ref{uniqueness} the map $\tau \co HX_p (X, \emptyset )\rightarrow HX_p' (X,\emptyset )$ is an isomorphism whenever the space $X$ is a finite coarse $CW$-complex.  The result now follows by looking at long exact sequences and applying the five lemma.
\end{proof}

\section{The coarse assembly map}

The ideas present in this article can be applied to the study of the assembly map present in the coarse Baum-Connes conjecture.  The coarse assembly map is described for proper metric spaces in \cite{HR2} and generalised in \cite{HPR}.  In this section we refine slightly the definition given in \cite{HPR} and prove that the refined coarse assembly map is an isomorphism for all finite coarse $CW$-complexes.  This result is an easy consequence of the machinery developed in the previous section together with some results presented in \cite{HPR}.

We begin by recalling some relevant definitions.

\begin{definition}
Let $X$ be a locally compact Hausdorff topological space, and let $A$ be a $C^\ast$-algebra.  Then an {\em $(X,A)$-module} is a Hilbert $A$-module $E$ equipped with a morphism $\varphi \co C_0(X)\rightarrow {\cal L}(E)$ of $C^\ast$-algebras.  The module $E$ is called {\em adequate} if $\overline{\varphi [C_0 (X)] E} = E$, and the operator $\varphi (f)$ is compact only when the function $f$ is the zero function.
\end{definition}

Here ${\cal L}({\cal E})$ denotes the $C^\ast$-algebra of all operators that admit adjoints on the Hilbert $A$-module $E$.  We usually omit explicit mention of the morphism $\varphi$ when talking about $(X,A)$-modules.  Note that when we talk about an operator between Hilbert $A$-modules being compact, we mean compact in the sense of operators between Hilbert modules, rather than compact as a bounded linear operator between Banach spaces.

\begin{definition}
Let $E$ be an $(X,A)$-module, and consider an operator $T\in {\cal L}(E)$.  Then we define the {\em support} of $T$, $\supp (T)$, to be the set of all pairs $(x,y)\in X\times X$ such that given functions $f\in C_0 (X)$ and $g\in C_0 (X)$, the equality $fTg = 0$ implies that either $f(x) = 0$ or $g(y) = 0$.
\end{definition}

\begin{definition}
Let $X$ be a coarse topological space, and let $E$ be an $(X,A)$-modules.  Consider an operator $T\in {\cal L}(E)$.

\begin{itemize}

\item The operator $T$ is said to be {\em locally compact} if the operators $Tf$ and $fT$ are both compact for all functions $f\in C_0(X)$.

\item The operator $T$ is said to be {\em pseudolocal} if the commutator $fT - Tf$ is compact for all functions $f\in C_0(X)$.

\item The operator $T$ is said to be {\em controlled} if the support, $\supp (T)$, is an entourage.

\end{itemize}
\end{definition}

\begin{definition} 
Let $X$ be a coarse topological space, and let $E$ be an adequate $(X,A)$-module.  Then we define the $C^\ast$-algebra $C^\ast_A (X)$ to be the $C^\ast$-algebra generated by all controlled and locally compact operators on the module $E$.  We define the $C^\ast$-algebra $D^\ast_A (X)$ to be the $C^\ast$-algebra generated by all controlled and pseudolocal operators on the module $E$.
\end{definition}

As our terminology suggests, the $(X,A)$-module $E$ which we are using in the above definition does not really matter to us.  The $K$-theory groups $K_n C^\ast_A (X)$ and $K_n D^\ast_A (X)$ do not depend on the precise choice of adequate $(X,A)$-module; see \cite{HPR} for a proof of this fact.

By proposition 5.5 of \cite{HPR} the assignments $X\mapsto K_n C^\ast_A (X)$ are coarsely-invariant functors on the category of coarse topological spaces.  
The $C^\ast$-algebra $C^\ast_A (X)$ is an ideal in the $C^\ast$-algebra $D^\ast_A (X)$.  We can thus form the quotient algebra $D^\ast_A (X)/C^\ast_A (X)$ and so we have an exact sequence:
$$\xymatrix@1{
K_{n+1} (D^\ast_A X) \ar[r] & K_{n+1}(D^\ast_A X / C^\ast_A X) \ar[r]^-{\alpha} & K_n (C^\ast_A X) \ar[r] & K_n(D^\ast_A X)
}$$

By proposition 7.3 of \cite{HPR} the $K$-theory group $K_{n+1}(D^\ast_A X / C^\ast_A X)$ is naturally isomorphic to the $KK$-theory group $KK^{-n} (C_0 (X),A)$.  The functor $X\mapsto KK^{-n} (C_0 (X),A)$ is a locally finite homology theory; it can be considered to be locally finite $K$-homology with coefficients in the $C^\ast$-algebra $A$.

\begin{proposition}
Let $X$ be a coarse paracompact topological space.  Let $KX_n (X;A)$ denote the coarse homology theory defined by coarsening the functor $X\mapsto KK^{-n}(C_0 (X),A)$.\footnote{See section \ref{chtsec}.}  Then there is a map $\alpha_\infty \co KX_n (X;A)\rightarrow K_n C^\ast_A (X)$ such that we have a commutative diagram:
$$\xymatrix{
{KK^{-n} (C_0 (X),A)} \ar[rd]^{\alpha} \ar[d]_{c} \\
{KX_n (X;A)} \ar[r]^{\alpha_\infty} & {K_n C^\ast_A (X)} \\
}$$
\end{proposition}

\begin{proof}
Let $({\cal U}_i )$ be an open coarsening family on the space $X$.  Recall from the comments at the end of section \ref{chtsec} that we can define a proper continuous map $\kappa \co X\rightarrow |{\cal U}_i |$ by the formula:
$$\kappa_i (x) = \sum_{U\in {\cal U}_i } \varphi_U (x) U$$

The coarsening map $c$ is defined to be the direct limit of the induced maps ${\kappa_i}_\ast \co KK^{-n}(C_0 (X),A)\rightarrow KK^{-n} ( C_0 (|{\cal U}_i |) , A)$.  The map $\alpha_\infty$ can be defined to be the direct limit of the maps:\footnote{The space $X$ is coarsely equivalent to each nerve $|{\cal U}_i |$, so the groups $K_n C^\ast_A (X)$ and $K_n C^\ast_A (|{\cal U}_i |)$ are isomorphic}
$$\alpha \co KK^{-n} (C_0 (|{\cal U}_i |) ,A)\rightarrow K_n C^\ast_A (X)$$
\end{proof}

The map $\alpha_\infty \co KX_n (X;A)\rightarrow K_n C^\ast_A (X)$ is called the {\em coarse assembly map}.  The {\em coarse Baum-Connes conjecture} asserts that this map is an isomorphism for all metric spaces of bounded geometry.\footnote{A metric space is said to have {\em bounded geometry} if it is coarsely equivalent to a discrete metric space $Y$ such that the supremum $\sup_{y\in Y} |B(y,r)|$ is finite for all real numbers $r>0$.  The coarse Baum-Connes conjecture is actually now known to be false in general; see \cite{HLS}.} 

The machinery we have developed in this article together with some results from \cite{HPR} tells us that the coarse assembly map is an isomorphism for all finite coarse $CW$-complexes.

To see this fact, first note that by an argument similar to that given to prove theorem 11.2 of \cite{HPR} the functors $K_n C^\ast_A$ are coarsely homotopy-invariant, and by corollary 9.5 of \cite{HPR} these functors satisfy the coarse excision axiom.  Hence the sequence of functors $\{ K_n C^\ast_A \}$ is a coarse homology theory.

We now need some computations.

\begin{lemma}
Let $+$ denote the one-point space.  Then the $K$-theory groups $K_n D^\ast_A (+)$ are all trivial.
\end{lemma}

\begin{proof}
The $C^\ast$-algebra $D^\ast_A (+)$ consists of all operators on some infinite-dim\-en\-sio\-nal Hilbert $A$-module.  An Eilenberg swindle argument tells us that the $K$-theory groups of such a $C^\ast$-algebra are all trivial.
\end{proof}

The map $\alpha \co KK^{-n}(C_0 (+),A)\rightarrow K_n C^\ast_A (+)$ fits into an exact sequence:
$$\xymatrix@1{
K_{n+1} D^\ast_A (+) \ar[r] & KK^{-n}(C_0 (+),A) \ar[r]^-{\alpha} & K_n C^\ast_A (+) \ar[r] & K_n D^\ast_A (+)
}$$

Further, the coarsening map $c\co KK^{-n} (C_0 (+) , A)\rightarrow KX_n (+;A)$ is clearly an isomorphism.  The above lemma therefore tells us that the coarse assembly map $\alpha_\infty \co KX_n (+;A)\rightarrow K_n C^\ast_A (+)$ is an isomorphism.

\begin{lemma}
Let $R$ be a generalised ray.  Then the groups $K_n C^\ast_A (R)$ are all trivial.
\end{lemma}

\begin{proof}
The generalised ray $R$ is flasque in the sense described in section 10 of \cite{HPR}.  Hence, by proposition 10.1 of \cite{HPR}, the $K$-theory groups $K_n C^\ast_A (R)$ are all trivial.
\end{proof}

\begin{lemma}
Let $R$ be a generalised ray.  Then the groups $KX_n (R;A)$ are all trivial.
\end{lemma}

\begin{proof}
We can find a coarsening family $({\cal U}_i )$ for the space $R$ such that each nerve $|{\cal U}_i |$ is properly homotopic to the ray $[0,\infty )$.  An Eilenberg swindle argument tells us that the groups $KK^{-n} (C_0 [0,\infty ) , A)$ are all trivial.  Taking direct limits, the groups $KX_n (R;A)$ must also be trivial.
\end{proof}

The above three lemmas tell us that the map $\alpha \co KK^{-n}(C(X);A)\rightarrow K_n C^\ast_A (X)$ is an isomorphism whenever the space $X$ is a single point or a generalised ray.  In particular, the coarse assembly map $\alpha_\infty$ is an isomorphism for such spaces.

\begin{theorem}
The coarse assembly map $\alpha_\infty \co KX_n (X;A)\rightarrow K_n C^\ast_A (X)$ is an isomorphism whenever the space $X$ is coarsely homotopy-equivalent to a finite coarse $CW$-complex.
\end{theorem}

\begin{proof}
The result is immediate from the previous three lemmas and theorem \ref{uniqueness}.
\end{proof}

A well-known descent argument (see \cite{Roe1}) enables us to deduce some results about the Novikov conjecture.

\begin{corollary}
Let $\Gamma$ be a finitely presented group such that the corresponding metric space $|\Gamma |$\footnote{The metric space $| \Gamma |$ is formed by equipping the group $\Gamma$ with the word-length metric.}, equipped with the bounded coarse structure, is coarsely homotopy-equivalent to a finite coarse $CW$-complex.  Then the group $\Gamma$ satisfies the Novikov conjecture.
\noproof
\end{corollary}

The above result is not new.  For example, any metric space coarsely homotopy-equivalent to a finite coarse $CW$-complex is certainly uniformly embeddable into a Hilbert space, so the above corollary is included in the main result of \cite{Yu2}.

However, another result is possible if we look at continuously controlled coarse structures rather than bounded coarse structures.  See for example \cite{HPR} and the final chapter of \cite{Roe1} for the descent arguments that are necessary in this case.

\begin{definition}
Let $X$ be a locally compact proper metric space, equipped with a compactification $\overline{X}$.  Let $X^b$ denote the space $X$ equipped with the bounded coarse structure, and let $X^\mathrm{cc}$ denote the space $X$ equipped with the continuously controlled coarse structure arising from the compactification $\overline{X}$.  Then we call the compactification $\overline{X}$ a {\em coarse compactification} if the identity map $1\co X^b \rightarrow X^\mathrm{cc}$ is a coarse map.
\end{definition}

\begin{corollary}
Let $\Gamma$ be a discrete group such that:

\begin{itemize}

\item The classifying space $B\Gamma$ can be represented as a finite complex.

\item The corresponding universal space $E\Gamma$ admits a coarse compactification, $\overline{E\Gamma}$, such that the space $E\Gamma^\mathrm{cc}$ is coarsely homotopy-equivalent to a finite coarse $CW$-complex, and the given $\Gamma$-action on the space $E\Gamma$ extends continuously to the compactification $\overline{E\Gamma}$.

\end{itemize}

Then the group $\Gamma$ satisfies the Novikov conjecture.
\noproof
\end{corollary}

Actually, the descent argument is sufficiently general to prove that the Novikov conjecture holds for any subgroup of a group of the kind described above.

\Addresses\recd

\end{document}